\documentclass[12pt,leqno]{article}

\usepackage{hyperref}
\usepackage{amsmath}
\usepackage{amssymb}
\usepackage{amsfonts}
\usepackage{amsopn}
\usepackage[all]{xy}
\usepackage{amsthm}
\usepackage{amscd}

\swapnumbers
\theoremstyle{plain}
\newtheorem{thm}{Theorem}[section]
\newtheorem{lem}[thm]{Lemma}

\theoremstyle{definition}
\newtheorem{ntt}[thm]{}
\newtheorem{rem}[thm]{Remark}
\newtheorem{dfn}[thm]{Definition}

\newcommand{\af}{\mathbb{A}}   
\newcommand{\pl}{\mathbb{P}}   
\newcommand{\zz}{\mathbb{Z}}   
\newcommand{\qq}{\mathbb{Q}}
\newcommand{\F}{{\mathrm{F}_4}}
\newcommand{\G}{\mathbb{G}}  

\newcommand{\id}{\mathrm{id}}       
\newcommand{\pr}{\mathrm{pr}}       
\newcommand{\rk}{\mathrm{rk}}

\newcommand{\Var}{\mathcal{V}ar} 
\newcommand{\M}{\mathrm{M}}     
\newcommand{\Mot}{\mathcal{M}}     
\newcommand{\Cor}{\mathcal{C}or} 

\DeclareMathOperator{\Spec}{\mathrm{Spec}}  
\DeclareMathOperator{\Pic}{\mathrm{Pic}}    
\DeclareMathOperator{\Gal}{\mathrm{Gal}}   
\DeclareMathOperator{\CH}{\mathrm{CH}}      
\DeclareMathOperator{\Ch}{\mathrm{Ch}}      
\DeclareMathOperator{\Mor}{\mathrm{Mor}}
\DeclareMathOperator{\End}{\mathrm{End}}
\DeclareMathOperator{\im}{\mathrm{Im}}      
\DeclareMathOperator{\SB}{\mathrm{SB}}      
\DeclareMathOperator{\MS}{\mathrm{MS}}
\DeclareMathOperator{\K}{\mathrm{K}}        
\DeclareMathOperator{\Nrd}{\mathrm{Nrd}}
\DeclareMathOperator{\GL}{\mathrm{GL}}
\DeclareMathOperator{\SL}{\mathrm{SL}}
\DeclareMathOperator{\diag}{\mathrm{diag}}
\DeclareMathOperator{\Gr}{\mathrm{Gr}}       
\DeclareMathOperator{\Seg}{\mathrm{Seg}}       
\DeclareMathOperator{\op}{\mathrm{op}}
\DeclareMathOperator{\res}{\mathrm{res}}
\DeclareMathOperator{\De}{\mathit{\Delta}}

\title{Motivic decomposition of a compactification of a Merkurjev-Suslin variety}

\author{N.~Semenov\footnote{
The author gratefully acknowledges the hospitality and support
of Bielefeld University.
Supported partially by CNRS, DAAD, and INTAS foundations.
}}

\date{}

\begin{document}

\maketitle

\begin{abstract}
We provide a motivic decomposition of a twisted form of a smooth
hyperplane section of $\Gr(3,6)$. This variety is a norm variety
corresponding to a symbol in $\K_3^M/3$. As an application we construct
a torsion element in the Chow group of this variety.

\end{abstract}

\section{Introduction}
In the present paper we study certain twisted forms of a smooth hyperplane
section of $\Gr(3,6)$. These twisted forms are smooth $\SL_1(A)$-equivariant
compactifications of a Merkurjev-Suslin variety corresponding to a central
simple algebra $A$ of degree $3$.
On the other hand,
these twisted forms are norm varieties corresponding
to symbols in $\K_3^M/3$ given by the Serre-Rost invariant $g_3$.
In the present paper we provide
a complete decomposition of the Chow motives of these varieties.

The history of this problem goes back to Rost and
Voevodsky. Namely, Rost obtained the celebrated decomposition of a norm quadric (see \cite{Ro98})
and later Voevodsky found some direct summand, called a generalized
Rost motive, in the motive of any norm variety (see \cite{Vo03}).
Note that the $\F$-varieties from \cite{NSZ05} can be considered
as a mod-$3$ analog of a Pfister quadric (more precisely, of
a maximal Pfister neighbour). In its turn, our variety can be
considered as a mod-$3$ analog of a norm quadric.

The paper is organized as follows. In section~\ref{notation}
we provide a background to the category of Chow motives. In section~\ref{decomposition}
we define 
a smooth compactification of a Merkurjev-Suslin variety $\MS(A,c)$ with
$A$ a central simple algebra of degree $3$,
describe its geometrical properties, and
decompose its Chow motive. 
Section~\ref{sectorsion} is devoted to an application of the obtained
motivic decomposition. Namely, using the ideas of Karpenko and Merkurjev
we construct a $3$-torsion element in the Chow group of our variety.

The main ingredients of our proofs are
results of Bia{\l}ynicki-Birula \cite{BB73}, Lefschetz hyperplane theorem,
and Segre embedding.

\section{Notation}\label{notation}
\begin{ntt}
The matrix notation of the present paper follows \cite{Inv}.

We use Galois descent language, i.e., identify a 
(quasi-projective) variety $X$ over a field $k$ with the variety $X_s=X\times_{\Spec k} \Spec k_s$ 
over the separable closure $k_s$ equipped with an action 
of the absolute Galois group 
$\Gamma=\Gal(k_s/k)$.
The set of $k$-rational points of $X$ is precisely the set of 
$k_s$-rational points of $X_s$ stable under 
the action of $\Gamma$.

We consider the Chow group $\CH^i(X)$ (resp. $\CH_i(X)$) of classes of algebraic cycles of 
codimension $i$ (resp. of dimension $i$) on an irreducible algebraic variety $X$ modulo rational equivalence
(see \cite{Ful}).

A Poincar\'e polynomial or generating function for a variety $X$ is, by definition, the polynomial
$\sum a_it^i\in\zz[t]$ with $a_i=\rk\CH^i(X)$.

The structure of the Chow ring of a Grassmann variety is of our
particular interest.
We do a lot of computations using formulae from Schubert calculus
(see \cite{Ful} 14.7).
\end{ntt}

Next we introduce the category of Chow motives over a field $k$ 
following \cite{Ma68} and \cite{CM04}. We remind the notion of 
a rational cycle and state the Rost Nilpotence
Theorem following \cite{CGM}.

\begin{ntt} Let $k$ be a field and $\Var_k$ be a category
of smooth projective varieties over $k$. Let $S$ denote any
commutative ring.
For any variety $X$ we set $\Ch(X):=\CH(X)\otimes_{\zz}S$.
First, we define the category of \emph{correspondences with
$S$-coefficients} (over $k$)
denoted by $\Cor_k(S)$. 
Its objects are smooth projective varieties over $k$. 
For morphisms, called correspondences,
we set $\Mor(X,Y):=\CH_{\dim X}(X\times Y)\otimes_{\zz}S$.
For any two correspondences $\alpha\in \Ch(X\times Y)$ and 
$\beta\in \Ch(Y\times Z)$ we define their composition 
$\beta\circ\alpha\in \Ch(X\times Z)$ as
\begin{equation}
\beta\circ\alpha ={\pr_{13}}_*(\pr_{12}^*(\alpha)\cdot \pr_{23}^*(\beta)),
\end{equation}
where $\pr_{ij}$ denotes the projection
on the $i$-th and $j$-th factors of $X\times Y\times Z$ 
respectively and ${\pr_{ij}}_*$, ${\pr_{ij}^*}$ denote
the induced push-forwards and pull-backs for Chow groups. 

The pseudo-abelian completion of $\Cor_k(S)$ is called the category
of \emph{Chow motives with $S$-coefficients} and is denoted by $\Mot_k(S)$.
The objects of $\Mot_k(S)$
are pairs $(X,p)$, where $X$ is a smooth projective variety
and $p\in\Mor(X,X)$ is an idempotent, that is, $p\circ p=p$. 
The morphisms between two objects $(X,p)$ and $(Y,q)$ 
are the compositions 
$q\circ \Mor(X,Y) \circ p$.
\end{ntt}

\begin{ntt}
By construction, $\Mot_k(S)$ is a tensor additive category,
where the tensor product is given by the usual product
$(X,p)\otimes(Y,q)=(X\times Y, p\times q)$.
For any cycle $\alpha$ we denote by $\alpha^t$ the
corresponding transposed cycle.
\end{ntt}

\begin{ntt} Observe that the composition product $\circ$
induces the ring structure on the abelian group $\Ch_{\dim X}(X\times X)$.
The unit element of this ring is the class of the diagonal map $\Delta_X$, which is defined by
$\Delta_X\circ \alpha = \alpha\circ \Delta_X=\alpha$ for all 
$\alpha\in \Ch_{\dim X}(X\times X)$.
The motive $(X,\Delta_X)$ will be denoted by $\Mot(X)$.
\end{ntt}

\begin{ntt} Consider the morphism
$(e,\id)\colon\{pt\}\times\pl^1\to\pl^1\times\pl^1$. 
Its image by means of
the induced push-forward $(e,\id)_*$ does not depend on the choice of the point
$e\colon\{pt\}\to \pl^1$
and defines the projector in $\CH_1(\pl^1\times\pl^1)$ denoted by $p_1$.
The motive $\zz(1)=(\pl^1,p_1)$ is called \emph{Lefschetz motive}.
For a motive $M$ and a nonnegative integer $i$ 
we denote its twist by $M(i)=M\otimes \zz(1)^{\otimes i}$.
\end{ntt}

\begin{ntt}\label{twistisom}
An isomorphism between the twisted motives $(X,p)(m)$ and $(Y,q)(l)$
is given by correspondences $j_1\in q\circ\Ch_{\dim X+l-m}(X\times Y)\circ p$ and
$j_2\in p\circ \Ch_{\dim Y+m-l}(Y\times X)\circ q$ such that
$j_1\circ j_2=q$ and $j_2\circ j_1=p$.
\end{ntt}

\begin{ntt}
Let $X$ be a smooth projective cellular variety.
The abelian group structure of $\CH(X)$
is well-known.
Namely, $X$ has a cellular filtration and 
the generators of Chow groups of
the bases of this filtration correspond to
the free additive generators of $\CH(X)$.
Note that the product of two cellular varieties
$X\times Y$ has a cellular filtration as well,
and $\CH^*(X\times Y)\cong \CH^*(X)\otimes \CH^*(Y)$
as graded rings.
The correspondence product of two cycles 
$\alpha=f_\alpha \times g_\alpha \in \Ch(X\times Y)$ and
$\beta=f_\beta \times g_\beta \in \Ch(Y\times X)$ is given
by (see \cite{Bo03} Lemma~5)
\begin{equation}
(f_\beta\times g_\beta)\circ(f_\alpha\times g_\alpha)=
\deg(g_\alpha \cdot f_\beta)(f_\alpha\times g_\beta),
\end{equation}
where $\deg\colon\Ch(Y)\to\Ch(\{pt\})=S$ is the degree map.
\end{ntt}

\begin{ntt}
Let $X$ be a projective variety of dimension $n$ over a field $k$. 
Let $k_s$ be the separable closure of $k$ and $X_s=X\times_{\Spec k} \Spec k_s$.
We say a cycle $J\in \Ch(X_s)$ is {\it rational}
if it lies in the image of the natural homomorphism
$\Ch(X)\to \Ch(X_s)$.
For instance, there is an obvious rational cycle $\Delta_{X_s}$ in 
$\Ch_n(X_s\times X_s)$ that is given by the diagonal class.
Clearly, all linear combinations, 
intersections and correspondence products of rational cycles
are rational.
\end{ntt}

\begin{ntt}[Rost Nilpotence]\label{nilpotence}
Finally, we shall also use the following fact (see \cite{CGM} Theorem~8.2)
called Rost Nilpotence theorem. 
Let $X$ be a projective homogeneous variety over $k$. Then
for any field extension $l/k$ the kernel of the natural
ring homomorphism $\End(\Mot(X))\to\End(\Mot(X_l))$ consists
of nilpotent elements.
\end{ntt}

\section{Motivic decomposition}\label{decomposition}
From now on we assume the characteristic of the base field $k$ is $0$.

It is well-known (see \cite{GH} Ch.~1, \S~5, p.~193) that the Grassmann variety
$\Gr(l,n)$ can be represented as the variety of $l\times n$ matrices
of rank $l$ modulo an obvious action of the group $\GL_l$.
Having this in mind we give the following definition.

\begin{dfn}
Let $A$ be a central simple algebra of degree $3$ over a field $k$, $c\in k^*$.
Fix an isomorphism $(A\oplus A)_s\simeq \M_{3,6}(k_s)$. Consider
the variety $D=D(A,c)$ obtained by Galois descent from the variety
$$
\{\alpha\oplus\beta\in (A\oplus A)_s\mid\rk(\alpha\oplus\beta)=3,\,
\Nrd(\alpha)=c\Nrd(\beta)\}/\GL_1(A_s),
$$
where $\GL_1(A_s)$ acts on $A_s\oplus A_s$ by the left multiplication.

This variety was first considered by M.~Rost.
\end{dfn}

Consider the Pl\"ucker embedding of $\Gr(3,6)$ into projective space
(see \cite{GH} Ch.~1, \S~5, p.~209). It is obvious that under this
embedding for all $c$ the variety $D(\M_3(k),c)$ is a
hyperplane section of $\Gr(3,6)$.

\begin{lem}
The variety $D$ is smooth.
\end{lem}
\begin{proof}(M.~Florence)
We can assume $k$ is separably closed.
Consider first the variety 
$$V=\{\alpha\oplus\beta\in\M_3(k)\oplus\M_3(k)=\M_{3,6}(k)\mid\rk(\alpha\oplus\beta)=3,\,
\det(\alpha)=c\det(\beta)\}.$$
An easy computation of differentials shows that $V$ is smooth.
The variety $V$ is a $\GL_3$-torsor over $D$ and, since $\GL_3$ is smooth,
this torsor is locally trivial for \'etale topology. Therefore
to prove its smoothness we can assume that this torsor is split.

Since $D\times_k\GL_3$ is smooth, $D\times_k\M_3$ is also smooth. Therefore
it suffices to prove that if $D\times_k\af^1$ is smooth, then $D$ is smooth.
But this is true for any variety.
Indeed, for any point $x$ on $D$ we have $T_{(x,0)}(D\times_k\af^1)=T_xD\oplus T_0\af^1=T_xD\oplus k$
and $\dim T_xD=\dim T_{(x,0)}(D\times_k\af^1)-1=\dim (D\times_k\af^1)-1=\dim D$.
\end{proof}

\begin{rem}
One can associate to the variety $D$ a Serre-Rost invariant $g_3(D)=[A]\cup [c]\in\mathrm{H}^3(k,\zz/3)$
(see \cite{Inv} \S~40).
This invariant is trivial if and only if $D$ is isotropic.
\end{rem}

It is easy to see that $D^0:=\MS(A,c):=\{a\in A\mid\Nrd(a)=c\}$ is an open orbit
under the natural right $\SL_1(A)$- or $\SL_1(A)\times\SL_1(A)$-action on $D$.
Namely, the open orbit
consists of all $\alpha\oplus\beta$ with $\rk(\alpha)=3$.
$D^0$ is called a Merkurjev-Suslin variety. In other words, the variety $D(A,c)$ is a smooth
$\SL_1(A)$-equivariant compactification of the Merkurjev-Suslin variety $\MS(A,c)$.

Denote as $\imath\colon D\to\SB_3(\M_2(A))$ the corresponding closed embedding.

\begin{lem}\label{lem3}
For the variety $D_s$ the following properties hold.

1. There exists a $\G_m$-action on $D_s$ with $18$ fixed points.
In particular, $D_s$ is a cellular variety.

2.
The generating function for $\CH(D_s)$ is equal to
$g=t^8+t^7+2t^6+3t^5+4t^4+3t^3+2t^2+t+1$. 

3. Picard group $\Pic(D_s)$ is rational.
\end{lem}
\begin{proof}
1. We can assume $c=1$.
The right action of $\G_m$ on $D_s$ is induced by the following action:
\begin{align*}
(\M_3(k_s)\oplus\M_3(k_s))\times\G_m&\to\M_3(k_s)\oplus\M_3(k_s)\\
(\alpha\oplus\beta,\lambda)&\mapsto
\alpha\diag(\lambda,\lambda^5,\lambda^6)\oplus
\beta\diag(\lambda^2,\lambda^3,\lambda^7)
\end{align*}
Note that this action is compatible with the left action of $\GL_3(k_s)$.

The $18$ fixed points of $D$ are the ${6\choose 3}=20$
$3$-dimensional standard subspaces of $\Gr(3,6)$ minus $2$ subspaces, generated by the first and
by the last $3$ basis vectors.

2. By the Lefschetz hyperplane theorem (see \cite{GH})
the pull-back $\imath_s^*$ is an isomorphism in
codimensions $i<\frac{\dim(\Gr(3,6))-1}2$.
Therefore $\rk\CH^i(D_s)=\rk\CH^i(\Gr(3,6))$ for such $i$'s.
Since Poincar\'e duality holds, we have
$\rk\CH_i(D_s)=\rk\CH_i(\Gr(3,6))$ for $i<\frac{\dim(\Gr(3,6))-1}2=4$.

It remains to determine
only the rank in the middle codimension. To do this observe that
$\rk\CH^*(D_s)=18$ (see \cite{BB73}). Therefore
$\rk\CH^4(D_s)=2\rk\CH^4(\Gr(3,6))-2=4$.

3. Consider the following commutative diagram:
\begin{equation}
\xymatrix{
\Pic(\SB_3(\M_2(A)))\ar[r]^-{\imath^*}\ar[d]&\Pic(D)\ar[d]^-{\res^*}\\
\Pic(\Gr(3,6))\ar[r]^-{\imath_s^*}&\Pic(D_s)}
\end{equation}
where the vertical arrows are the morphisms of scalar extension.
By the Lefschetz hyperplane theorem the map $\imath_s^*$ restricted
to $\Pic(\Gr(3,6))$ is an isomorphism.
Since $\Pic(\SB_3(\M_2(A)))$ is rational (see \cite{MT95}  and
\cite{NSZ05} Lemma~4.3), i.e., the
left vertical arrow is an isomorphism, the restriction map $\res^*$ is surjective.
On the other hand, it follows from a Hochschild-Serre spectral sequence
(see \cite{Ar82} \S~2) that $\Pic(D)$ can be identified with a subgroup
of $\zz$. We are done.
\end{proof}

\begin{rem}
It immediately follows from this Lemma that the variety $D$ is not a twisted
flag variety. Indeed, the generating functions of all twisted flag varieties over
a separabely closed field are
well-known and all of them are different from the generating function of $D_s$.
\end{rem}

\begin{ntt}\label{multrule}
We must determine a (partial) multiplicative structure of $\CH(D_s)$.
By Lefschetz hyperplane theorem the generators in
codimensions $0$, $1$, $2$, and $3$ are pull-backs
of the canonical generators $\De_{(0,0,0)}$, $\De_{(1,0,0)}$,
$\De_{(1,1,0)}$, $\De_{(2,0,0)}$, $\De_{(1,1,1)}$, $\De_{(2,1,0)}$,
$\De_{(3,0,0)}$ of $\Gr(3,6)$ (see \cite{Ful} 14.7). 
We denote these pull-backs as $1$, $h_1$, $h_2^{(1)}$, $h_2^{(2)}$,
$h_3^{(1)}$, $h_3^{(2)}$, and $h_3^{(3)}$ respectively.
In the codimension $4$ the pull-back is injective and the pull-backs 
$h_4^{(1)}:=\imath_s^*(\De_{(2,1,1)})$, $h_4^{(2)}:=\imath_s^*(\De_{(2,2,0)})$,
$h_4^{(3)}:=\imath_s^*(\De_{(3,1,0)})$, where $\imath$
is as above, form a subbasis of $\CH^4(D_s)$.

Consider the following diagram:
$$
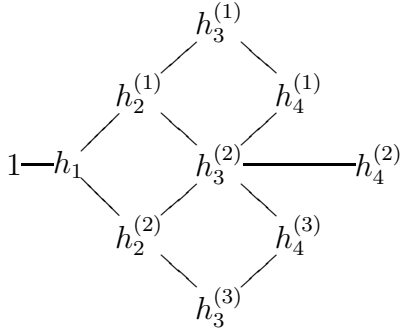
\xymatrix@M=0pt@=1em{
 & & & h_3^{(1)}\ar@{-}[rd] &  & &\\
 & &h_2^{(1)}\ar@{-}[ru] \ar@{-}[rd] & &h_4^{(1)} & &\\
1\ar@{-}[r] & h_1 \ar@{-}[ru] \ar@{-}[rd]& &h_3^{(2)}\ar@{-}[ru] \ar@{-}[rd]\ar@{-}[rr]& &h_4^{(2)}&\\
 & &h_2^{(2)}\ar@{-}[ru] \ar@{-}[rd]& & h_4^{(3)}& &\\
 & & & h_3^{(3)}\ar@{-}[ru] &  &  &
}
$$

Since pull-backs are ring homomorphisms, it immediately follows that
$$
h_1\cdot u = \sum_{u\to v} v,
$$ 
where $u$ is a vertex on the diagram, which
corresponds to a generator of codimension less than $4$, and 
the sum runs through all the edges going from $u$ 
one step to the right.
\end{ntt}

Next we compute some products in the middle codimension.

Since $\De_{(3,1,0)}\De_{(2,1,1)}=\De_{(2,2,0)}^2=0$ and
$\De_{(2,1,1)}^2=\De_{(3,1,0)}^2=\De_{(2,2,0)}\De_{(2,1,1)}=\De_{(2,2,0)}\De_{(3,1,0)}=\De_{(3,3,2)}$
(see \cite{Ful} 14.7),
we have $h_4^{(1)}h_4^{(3)}=(h_4^{(2)})^2=0$ and
$(h_4^{(1)})^2=(h_4^{(3)})^2=h_4^{(2)}h_4^{(3)}=h_4^{(1)}h_4^{(2)}=\imath_s^*(\De_{(3,3,2)})=pt$,
where $pt$ denotes the class of a rational point on $D_s$.

The next theorem shows that the Chow motive of $D$ with $\zz/3$-coefficients
is decomposable. Note that for any cycle $h$ in $\CH(D_s)$ or in $\CH(D_s\times D_s)$
the cycle $3h$ is rational.
\begin{thm}\label{mainthm}
Let $A$ denote a central simple algebra of degree $3$ over
a field $k$, $c\in k^*$, and $D=D(A,c)$.
Then
$$
\Mot(D)\simeq R\oplus(\oplus_{i=1}^5 R'(i)),
$$
where $R$ is a motive such that
over a separably closed field it becomes isomorphic
to $\zz\oplus\zz(4)\oplus\zz(8)$
and $R'\simeq\Mot(\SB(A))$.
\end{thm}
\begin{proof}
Consider the following commutative diagram (see \cite{CPSZ05} 5.5):
\begin{equation}
\footnotesize
\xymatrix{
D_s\times\pl^2\ar[r]^-{\imath_s\times\id_s}\ar[d]&\Gr(3,6)\times\pl^2\ar[r]^-{\Seg_s}\ar[d]&\Gr(3,18)\ar[d]\\
D\times\SB(A^{\op})\ar[r]^-{\imath\times\id}&\SB_3(\M_2(A))\times\SB(A^{\op})\ar[r]^-{\Seg}&\SB_3(\M_2(A)\otimes_k A^{\op})
}
\end{equation}
where the right horizontal arrows are Segre embeddings given
by the tensor product of ideals (resp. linear subspaces) and the vertical
arrows are canonical maps induced by the scalar extension $k_s/k$.

This diagram induces the commutative diagram of rings
\begin{equation}
\footnotesize
\xymatrix{
\Ch^*(D_s\times\pl^2)&\Ch^*(\Gr(3,6)\times\pl^2)\ar[l]^-{(\imath_s\times\id_s)^*}&\Ch^*(\Gr(3,18))\ar[l]^-{\Seg_s^*}\\
\Ch^*(D\times\SB(A^{\op}))\ar[u]&\Ch^*(\SB_3(\M_2(A))\times\SB(A^{\op}))\ar[l]^-{(\imath\times\id)^*}\ar[u]&
\Ch^*(\SB_3(\M_2(A)\otimes_k A^{\op}))\ar[u]^-{\simeq}\ar[l]^-{\Seg^*}
}
\end{equation}
Observe that the right vertical arrow is an isomorphism
since $\M_2(A)\otimes_k A^{\op}$ splits.

Let $\tau_3$ and $\tau_1$ be tautological vector bundles on $\Gr(3,6)$ and
$\pl^2$ respectively and let $e$ denote the Euler class (the top Chern class).
By \cite{CPSZ05} Lemma~5.7 the cycle
$(\imath_s\times\id_s)^*(e(\pr_1^*\tau_3\otimes\pr_2^*\tau_1))\in\Ch(D_s\times\pl^2)$
is rational.
A straightforward computation (cf. \cite{CPSZ05} 5.10 and 5.11) shows that
$r:=-(\imath_s\times\id_s)^*(e(\pr_1^*\tau_3\otimes\pr_2^*\tau_1))=
h_3^{(1)}\times 1+h_2^{(1)}\times H+h_1\times H^2\in\Ch^3(D_s\times\pl^2)$, where
$H$ is the class of a smooth hyperplane section of $\pl^2$.

Define following rational cycles
$\rho_i=r(h_1^i\times 1)\in\Ch^{3+i}(D_s\times\pl^2)$ for $i=1,\ldots,4$,
$\rho_0=r+h_1^3\times 1\in\Ch^3(D_s\times\pl^2)$ and
$\rho'_1=r(h_1\times 1)+h_1^4\times 1$.
A straightforward computation using multiplication rules~\ref{multrule}
shows that 
$(-\rho'_1)\circ\rho_3^t$ as well as
$(-\rho_{4-i})\circ\rho_i^t\in\Ch_2(\pl^2\times\pl^2)$ is the diagonal $\Delta_{\pl^2}$.
Moreover, the opposite compositions $(-\rho_0)^t\circ\rho_4$,
$(-\rho_1)^t\circ\rho_3$, $(-\rho_2)^t\circ\rho_2$,
$(-\rho_3)^t\circ\rho'_1$, and $(-\rho_4)^t\circ\rho_0$ give rational
pairwise orthogonal idempotents in $\CH_8(D_s\times D_s)$.

To finish the proof of the theorem it remains by \ref{twistisom} to lift
all these rational cycles $\rho_i$, $\rho_j^t$ to
$\Ch(D\times\SB(A^{\op}))$ and to
$\Ch(\SB(A^{\op})\times D)$ respectively in such a way that the corresponding
compositions of their preimages would give the diagonal $\Delta_{\SB(A^{\op})}$.

Fix an $i=0,\ldots,4$. Consider first any preimage
$\alpha\in\Ch(D\times\SB(A^{\op}))$ of $-\rho_{4-i}$ and
any preimage $\beta\in\Ch(\SB(A^{\op})\times D)$ of $\rho_i^t$.
The image of the composition $\alpha\circ\beta$ under the restriction map
is the diagonal $\Delta_{\pl^2}$. Therefore by Rost Nilpotence theorem
for Severi-Brauer varieties (see \ref{nilpotence}) $\alpha\circ\beta=\Delta_{\SB(A^{\op})}+n$, where
$n$ is a nilpotent element in $\End(\Mot(\SB(A^{\op})))$.
Since $n$ is nilpotent $\alpha\circ\beta$ is invertible and
$((\Delta_{\SB(A^{\op})}+n)^{-1}\circ\alpha)\circ\beta=\Delta_{\SB(A^{\op})}$.
In other words,
we can take 
$(\Delta_{\SB(A^{\op})}+n)^{-1}\circ\alpha$
as a preimage of $-\rho_{4-i}$ and
$\beta$
as a preimage of $\rho_i^t$. Note that $n$ is always a torsion element
and since $\End(\Mot(\SB(A^{\op})))\simeq\Mor(\Mot(\SB(A)),\Mot(\SB(A^{\op})))$
and $\CH(\SB(A))$ has no torsion, projective bundle theorem implies that
in fact $n=0$.

Denote as $R$ the remaining direct summand of the motive of $D$.
Comparing the left and the right hand sides of the decomposition over $k_s$
it is easy to see
that $R_s\simeq\zz\oplus\zz(4)\oplus\zz(8)$.
\end{proof}

\begin{rem}
Using a bit more messy computations one can show that
the same proof works for the motive of $D$ with integral coefficients.
\end{rem}

\section{Torsion}\label{sectorsion}
In this section we use Steenrod operations modulo $3$
(see \cite{Br03}, \cite{KM02} \S~3, and \cite{Me03}).
We denote the total Steenrod operation by $S^\bullet=S^0+S^1+\ldots$.

Let $X$ be a smooth projective variety over $k$.
For any cycle $p\in\CH(X\times X)$ we define its realization
$p_\star\colon\CH(X)\to\CH(X)$ as $p_\star(\alpha)=\pr_{2*}(\pr_1^*(\alpha)p)$,
$\alpha\in\CH(X)$, where $\pr_1,\pr_2\colon X\times X\to X$
denote the first and the second projections.
As $\deg\colon\CH_0(X)\to\zz$ we denote the usual degree map.

The goal of the present section is to prove the following theorem.
\begin{thm}\label{torsion}
Assume that the variety $D$ is anisotropic. Then $\CH_2(D)$ contains
$3$-torsion.
\end{thm}

\begin{ntt}
The proof of this Theorem consists of several parts. First we define
an important element $d$ as follows. 
The kernel of the push-forward map
$$(\imath_s)_*|_{\CH_4(D_s)}\colon\CH_4(D_s)\to\CH_4(\Gr(3,6))$$ has rank $1$, since
by Lefschetz hyperplane theorem the push-forward
$((\imath_s)_*|_{\CH_4(D_s)})\otimes\qq$
is surjective. Denote as $d\in\CH_4(D_s)$ a generator of this kernel.
Projection formula immediately implies that $(\imath_s)_*(\alpha d)=0$
for any $\alpha\in\im\imath_s^*$ and therefore
by Lefschetz hyperplane theorem $\alpha d=0$.

From now on we work with Chow groups modulo $3$.
\end{ntt}

\begin{lem}\label{lem1}
We have 
\begin{enumerate}
\item $d^2\ne 0\mod 3$,
\item the total Chern class of the tangent bundle 
$$c(-T_{D_s})=1+h_1+h_1^2-h_1^3-h_1^4-h_1^5,$$ and
\item $S^\bullet(d)=d$.
\end{enumerate}
\end{lem}
\begin{proof}
The first equality is just a routine computation, which uses
Poincar\'e duality on $\CH(D_s)$.

Next we compute the total Chern class of the tangent bundle $T_{D_s}$.
Since $D_s$ is a hyperplane section of $\Gr(3,6)$ we have the following
exact sequence:
$$0\to T_{D_s}\to\imath_s^*(T_{\Gr(3,6)})\to\imath_s^*(\mathcal{O}_{\Gr(3,6)}(1))\to 0$$
Therefore $c(T_{D_s})\imath_s^*(c(\mathcal{O}_{\Gr(3,6)}(1)))=\imath_s^*(c(T_{\Gr(3,6)}))$.
Since $\imath_s^*(c(\mathcal{O}_{\Gr(3,6)}(1)))=1+h_1$ and
$\imath_s^*(c(T_{\Gr(3,6)}))=1-h_1^2-h_1^3+h_1^5$, we have $c(T_{D_s})=1-h_1-h_1^3+h_1^4$
and $c(-T_{D_s})=1+h_1+h_1^2-h_1^3-h_1^4-h_1^5$.

To prove the last assertion note that
$\Delta_{D_s}=\Delta'\pm d\times d$, where $\Delta'$ is
a part of the diagonal $\Delta_{D_s}$, which does not involve $d$, i.e.,
which comes from $\Gr(3,6)$. Let $\delta\colon D_s\to D_s\times D_s$
denote the diagonal morphism.

Now
$$S^\bullet(\pm d\times d)=S^\bullet(\Delta_{D_s}-\Delta')=
S^\bullet(\delta_*(1)-\Delta')=S^\bullet(\delta_*(1))-S^\bullet(\Delta').$$
To prove that $S^\bullet(d)=d$ we must show that the right hand side
does not contain summands of the form $d\times\alpha$, $\alpha\in\Ch(D_s)$,
different from $\pm d\times d$. Therefore the summand $S^\bullet(\Delta')$
is not interesting for us.

We have
\begin{align*}
S^\bullet(\delta_*(1))=c(T_{D_s\times D_s})\delta_*(S^\bullet_{D_s}(1)c(-T_{D_s}))
=c(T_{D_s\times D_s})(c(-T_{D_s})\times 1)\delta_*(1)\\
=c(T_{D_s\times D_s})(c(-T_{D_s})\times 1)\Delta_{D_s},
\end{align*}
where the second equality follows from projection formula.
But by item~2. the Chern classes $c_i(T_{D_s})$
don't involve $d$, i.e., lie in the image of $\imath_s^*$.
The lemma is proved.
\end{proof}

In the notation of Theorem~\ref{mainthm} denote as $p\in\Ch_8(D\times D)$ the
projector corresponding to the motive $R$, i.e., $R=(D,p)$.
From the proof of Theorem~\ref{mainthm} it is easy to see that
$$p_s=1\times pt\pm d\times d+pt\times 1.$$
Since the natural map $\Pic(D)\to\Pic(D_s)$ is an isomorphism
(see Lemma~\ref{lem3}(3)), we denote as $h_1$ the canonical generator
of $\Pic(D_s)$ as well as the correponding generator of $\Pic(D)$.

\begin{lem}\label{lem2}
The following properties of $D$ hold:
\begin{enumerate}
\item The natural group homomorphism $\CH_0(D)\to\CH_0(D_s)$ is injective.
Its image is generated by zero cycles of degree divisible by $3$.
\item $S^1(p_\star(h_1^6))=h_1^8$.
\end{enumerate}
\end{lem}
\begin{proof}

1. By \cite{CM06} Theorem~6.5 it suffices to show that the class $\mathcal{A}(D)$ of
all field extensions $E/k$ such that $D(E)\ne\emptyset$ is connected and for any $L\in\mathcal{A}(D)$ the
group $\CH_0(D_L)=\zz$. The first assertion follows from
\cite{CM06} Theorem~11.3, since any field extension $E/k$ such that $D=D(A,c)$
has an $E$-point splits the Jordan algebra $J(A,c)$ obtained by the first Tits construction,
and vice versa.

To prove that $\CH_0(D_L)=\zz$ for any $L\in\mathcal{A}(D)$
it suffices to check that for any field extension
$E/L$ any two rational points
of $D_E$ are rationally equivalent (see \cite{CM06} Lemma~5.2).
If the algebra $A_E$ is not split, then all rational points of $D_E$
are contained in $\MS(A_E,1)\simeq\SL_1(A_E)$. Since $\SL_1(A)$ is rational and
homogeneous, this implies that $D_E$ is $R$-trivial, and, hence,
$\CH_0(D_E)=\zz$. If the algebra $A_E$ splits, then obviously $\CH_0(D_E)=\zz$.

2. 
The proof of this item is similar to the proof of Corollary~4.9 \cite{KM02}.
By \cite{KM02} Lemma~3.1
\begin{equation}\label{form1}
S^\bullet(p_\star(h_1^6))=S^\bullet(p)_\star(h_1^6(1+h_1^2)^6c(-T_D)).
\end{equation}
Therefore $S^1(p_\star(h_1^6))$ equals the $0$-dimensional component of the
right hand side. Assume that 
\begin{equation}\label{form2}
S^1(p)_\star(h_1^6)=0.
\end{equation}
Then an easy computation using item~1. shows that the right hand side
of (\ref{form1}) equals $p_\star(h_1^8)=h_1^8$.

To prove (\ref{form2}) it suffices to show that $\deg S^1(p)_\star(h_1^6)$ is
divisible by $9$ (cf. \cite{KM02} Proof of Corollary~4.5).
Without loss of generality we can compute this degree over $k_s$.
It follows that $\deg S^1(p)_\star(h_1^6)=\deg h_1^6\pr_{1*}(S^1(p_s))$
(see \cite{KM02} Proof of Corollary~4.5). But $\pr_{1*}(S^1(p_s))$ is divisible
by $3$ (see Lemma~\ref{lem1}(3)) and for any $\alpha\in\Ch_2(D_s)$ the
product $h_1^6\alpha$ is divisible by $3$. We are done.
\end{proof}

Now we are ready to prove Theorem~\ref{torsion}.
Consider the cycle $S^1(p_\star(h_1^6))$. Since $\deg h_1^8=42$ and $D$ is
anisotropic, by lemma~\ref{lem2} this cycle is non-zero.
Therefore $p_\star(h_1^6)\in\Ch_2(D)$ is non-zero.
On the other hand, $(p_s)_\star(h_1^6)=0$. In other words, $p_\star(h_1^6)$
is a non-trivial torsion element in $\Ch_2(D)$.

\subparagraph{Acknowledgements\\}
I express my sincerely gratitude to M.~Florence, N.~Karpenko, and I.~Panin for interesting
discussions concerning the subject of the present paper.

\noindent
N.~Semenov\\
Fakult\"at f\"ur Mathematik\\
Universit\"at Bielefeld\\
Deutschland

\end{document}